\documentclass[draft]{imsart}

\usepackage{amsthm,amsmath,amssymb,natbib}
\RequirePackage[dvips]{hyperref}

\startlocaldefs

\theoremstyle{plain}
\newtheorem{theorem}[equation]{Theorem}
\newtheorem{lemma}[equation]{Lemma}
\newtheorem{corollary}[equation]{Corollary}

\theoremstyle{definition}
\newtheorem{definition}[equation]{Definition}
\newtheorem{assumption}{Assumption}

\numberwithin{equation}{section}

\theoremstyle{remark}
\newtheorem{remark}[equation]{Remark}

\newcommand{\F}{\mathcal{F}}

\newcommand{\T}{\mathcal{T}}

\newcommand{\X}[1]{\mathcal{X}_{#1}(x_1^n)}
\newcommand{\V}[1]{V_{#1}(x_1^n)}

\newcommand{\N}{\mathbb{N}}

\newcommand{\Z}{\mathbb{Z}}

\renewcommand{\P}{\mathbb{P}}
\newcommand{\Id}{{\mathbf {1}}}
\newcommand{\given}{|}
\newcommand{\ee}{e^{\frac1e}}

\newcommand{\ML}{{\textup{ML}}}

\newcommand{\Int}{{\textup{Int}}}
\newcommand{\argmin}[1]{\underset{#1}{\textup{arg min}}}

\newcommand{\sft}{\F^{d}(x_1^n)}
\renewcommand{\tau}{\T}

\endlocaldefs

\begin{document}

\begin{frontmatter}

\title{Some upper bounds for the rate of convergence of penalized likelihood context tree estimators}

\runtitle{Upper bounds for PL context tree estimators}

\author{\fnms{Florencia} \snm{Leonardi}\corref{}\ead[label=e1]
{florencia@usp.br}\thanksref{t1}}
\thankstext{t1}{This work is part of PRONEX/FAPESP's
project \emph{Stochastic behavior, critical phenomena and rhythmic pattern
identification in natural languages} (grant number 03/09930-9), 
CNRS-FAPESP project \emph{Probabilistic phonology of rhythm} and
CNPq project \emph{Rhythmic patterns, prosodic domains and probabilistic modeling
in Portuguese Corpora} (grant number 485999/2007-2).
} 
\address{Instituto de Matem\'atica e Estat\'istica, Universidade
de S\~ao Paulo\\Rua do Mat\~ao 1010 - Cidade Universit\'aria\\
CEP 05508-090 - S\~ao Paulo - SP - Brazil\\ \printead{e1}}
\affiliation{Instituto de Matem\'atica e Estat\'istica, Universidade
de S\~ao Paulo}

\runauthor{Florencia Leonardi}  
  
\begin{abstract}
We find upper bounds for the probability of underestimation and overestimation
errors in penalized likelihood context tree estimation.  
The bounds are explicit and applies to processes of not necessarily finite memory.
We allow for general penalizing terms and we give conditions over 
the maximal depth of the estimated trees in order to get strongly consistent estimates. 
This generalizes previous results obtained in the case of estimation of the order of a Markov chain.
\end{abstract}



\end{frontmatter}

\section{Introduction}

In this paper we obtain an exponential upper bound for the underestimation of the
\emph{context tree} of a variable memory process by penalized likelihood (PL)
criteria and a sub-exponential upper bound for the overestimation event. 
Our result applies to processes of not necessarily finite memory that satisfies some 
continuity requirements, generalizing the bound obtained in \citet{dorea2006} for the
estimation of the order of a Markov chain by similar methods (EDC
criterion).  

The concept of context tree was first introduced by
\citet{rissanen1983} to denote the minimum set of sequences that are necessary
to predict the next symbol in a finite memory stochastic chain.  A particular
case of context tree is the set of all sequences of length $k$, representing a
Markov chain of order $k$. For that reason, context trees allow a more
detailed and parsimonious representation of processes than finite order Markov
chains do.

In the statistical literature, the processes allowing a context tree
representation are called Variable Length Markov Chains 
\citep{buhlmann1999}. This class of models has shown to be useful in real
data modeling, as for example, for the case of protein classification into families
\citep{bejerano2001a, leonardi2006a}.

Historically, the estimation of the context tree of a process has
been addressed by different versions of the algorithm \emph{Context},
introduced by Rissanen in its seminal paper. This
algorithm was proven to be weak consistent in the case of bounded
memory  \citep{buhlmann1999} and also in the case of unbounded memory 
\citep{ferrari2003,duarte2006}.  Recently, in \citet{galves2006} 
it was obtained an upper bound for the rate of convergence of the
algorithm Context in the case of
bounded memory processes. A generalization of this result to the case of
unbounded memory processes was given in \citet{galves2008}. 

The estimation of context trees by PL criteria 
had not been addressed in the literature until the recent work by
\citet{csiszar2006}.  The reason for that was the exponential cost of the
estimation, due to the number of trees that had to be considered in order
to find the optimal one. In their article, Csisz\'ar and Talata 
showed that the Bayesian Information Criterion
(BIC), which is a particular case of the PL estimators (using a penalizing
term growing logarithmically), is strongly consistent and can be computed in 
linear time, using a suitable version of the Context Tree Weighting method
of Willems, Shtarkov and Tjalkens \citep{willems1995,willems1998}. 
Their result applies to unbounded memory processes
and the depth of the estimated tree is allowed to grow with the sample size as
a sub-logarithmic function. 
This last condition was proven to be unnecessary
in the case of finite memory processes, as proven in \citet{garivier2006}.
An explicit bound on the rate of convergence of the PL context tree estimators had remained until 
now as an open question. 

The paper is organized as follows. In Section~2 we introduce some
definitions and state the main result. In Section~3 we present the proofs and in Section~4
we do some final remarks. Finally, Section~5 constitutes and appendix that contains some
results needed in our proofs and obtained elsewhere in the literature.

\section{Definitions and results}

In what follows $A$ will represent a finite alphabet of size $|A|$.
Given two integers $m\leq n$, we will denote by $w_m^n$ the sequence
$(w_m, \ldots, w_n)$ of symbols in $A$. The length of the sequence
$w_m^n$ is denoted by $\ell(w_m^n)$ and is defined by $\ell(w_m^n) =
n-m+1$.  Any sequence $w_m^n$ with $m > n$ represents the empty string
and is denoted by $\lambda$. The length of the empty string is
$\ell(\lambda) = 0$. In the sequel $A^j$ will denote the set of all sequences of length $j$
over $A$.

Given two sequences $w = w_m^n$ and $v = v_j^k$, we will denote by
$vw$ the sequence of length $\ell(v) + \ell(w) $ obtained by
concatenating the two strings.  In particular, $\lambda w = w\lambda=
w$.  The concatenation of sequences is also extended to the case in
which $v$ denotes a semi-infinite sequence, that is
$v = (\dotsc, v_{-2},v_{-1})$, denoted by $v=v_{-\infty}^{-1}$.

We say that the sequence $s$ is a \emph{suffix} of the sequence $w$ if
there exists a sequence $u$, with $\ell(u)\geq 1$, such that $w = us$.
In this case we write $s\prec w$. When $s\prec w$ or $s=w$ we write
$s\preceq w$.

\begin{definition} 
A set  $\tau$ of finite or semi-infinite sequences is a \emph{tree} if no sequence
$s \in \tau$ is a suffix of another sequence $w \in \tau$. This
property is called the \emph{suffix property}.
\end{definition}

We define the \emph{height} of the tree $\tau$ as
\begin{equation*}
h(\tau) = \sup\{\ell(w) : w\in\tau\}. 
\end{equation*}

In the case $h(\tau)<+\infty$ 
we say that $\tau$ is \emph{bounded}
and we denote by $|\tau|$ the number of sequences in $\tau$. 
On the other hand, if $h(\tau)=+\infty$ 
we say that the tree $\tau$ is
\emph{unbounded}.
 
Given a tree $\tau$ and an integer $K$ we will denote by $\tau|_K$ the
tree $\tau$ \emph{truncated} to level $K$, that is
\begin{equation*}
\tau|_K = \{w \in \tau\colon \ell(w) \le K\} \cup \{ w\colon \ell(w)=K 
\text{ and } w\prec u, \text{ for some } u\in\tau\}.
\end{equation*}
The expression $\Int(\tau)$ will denote the set of all sequences
that are suffixes of some $u\in\tau$, that is
\begin{equation*}
\Int(\tau) = \{w\colon w\prec u, \text{ for some }u\in\tau\}.
\end{equation*}

We will say that a tree $\tau$ is \emph{complete} if  for every semi-infinite sequence
$w_{\infty}^{-1}$ there exists a sequence $s\in\tau$ such that $s\preceq w_{-\infty}^{-1}$.

Consider a stationary ergodic stochastic chain $\{X_t: t\in\Z\}$
over $A$. Given a sequence $w\in A^j$ we denote by 
\begin{equation*}
 p(w) \,=\,  \P(X_1^j = w)
\end{equation*}
the stationary probability of the cylinder defined by the sequence $w$.
If $p(w) > 0$ we write
\begin{equation*}
p(a|w) \,=\, \P ( X_0 =a \given X_{-j}^{-1}=w)\,.
\end{equation*}

In the sequel we will use the simpler notation $X_t$ for the process 
$\{X_t: t\in\Z\}$.

\begin{definition}\label{context}
A sequence $w\in A^j$ is a \emph{context} for the process $X_t$ if it satisfies 
\begin{enumerate}
\item For any semi-infinite sequence $x_{-\infty}^{-1}$ having $w$ 
as a suffix  
\begin{equation*}
\P ( X_0 =a \given X_{-\infty}^{-1}=x_{-\infty}^{-1}) \,=\, p(a|w),
\quad\text{for all $a\in A$}.
\end{equation*}
\item No suffix of $w$ satisfies (1).
\end{enumerate}
An \emph{infinite context} is a semi-infinite
sequence $w_{-\infty}^{-1}$ such that any of its suffixes $w_{-j}^{-1}$, $j=1,2,\dotsc$  is a context. 
\end{definition}

Definition~\ref{context} implies that the set of all contexts (finite or
infinite) satisfies the suffix property and hence it is a tree.  This tree is called
\emph{the context tree} of the process $X_t$ and will be denoted by
$\tau_0$. 

\begin{remark}
In this paper we will also consider i.i.d. processes. We will assume
that these processes are compatible with a particular tree, given by
the set $\{\lambda\}$.
\end{remark}

Define the sequence $\{\alpha_k\}_{k\in\N}$ as
\begin{align}\label{alpha}
\alpha_0 &:=  \inf_{w\in\tau_0,a\in A} \{\, p(a|w) \,\},\notag \\
\alpha_k &:= \inf_{u\in A^k}\;\sum_{a\in A}\;\inf_{w\in\tau_0, w\succ u}\{\,p(a|w)\,\}.
\end{align}

\begin{assumption}
From now on we will assume the process $X_t$ satisfies 
\begin{enumerate}
\item $\alpha_0 > 0$ and
\item $\alpha := \sum_{k\in\N} (1-\alpha_k) \;< \;+\infty$.
\end{enumerate}
\end{assumption}

The positivity assumption over $\alpha_0$ implies that the context tree of the
process $X_t$ is complete, i.e., any semi-infinite sequence $w_{-\infty}^{-1}$
belongs to $\tau_0$ or has a suffix that belongs to $\tau_0$.  
The second assumption is related to the loss of memory of a process of infinite order. 
(see \citet{galves2008} for more details). 

In what follows we will assume $x_1, x_2, \dotsc, x_n$ is a
sample of the process $X_t$. Let $d(n) < n$ be a function taking integer values and growing to
infinity with $n$. 
This will denote the maximal height of
the estimated context trees (and will be denoted simply by $d$). Then, given a sequence $w$, with $1\leq
\ell(w) \leq d$, and a symbol $a\in A$ we denote by $N_n(w,a)$ the
number of occurrences of symbol $a$ preceded by the sequence $w$,
starting at $d+1$, that is,
\begin{equation} \label{eq:Nn}
N_n(w,a)=\sum_{t=d+1}^{n}\Id\{x_{t-\ell(w)}^{t-1}=w,x_t = a\}. 
\end{equation}
On the other hand, $N_n(w)$ will denote the sum $\sum_{a \in A}
N_n(w,a)$.

\begin{definition}\label{feasible}
We will say that the tree $\tau$ is \emph{feasible} if it is complete,  $h(\tau)\leq
d$, $N_n(w)\geq 1$ for all $w\in\tau$ and any string $w'$ with
$N_n(w')\geq 1$ either belongs to $\tau$, is a suffix of some $w\in\tau$ or has a
suffix $w$ that belongs to $\tau$.
\end{definition}

We will denote by $\sft$ the set of all feasible trees. Then, given a
tree $\tau\in\sft$, the maximum likelihood of the sequence
$x_1,\dotsc,x_n$ is given by
\begin{equation}\label{pmle}
\hat\P_{\ML,\tau}(x_1^n) =  \prod_{w\in \tau}\prod_{a\in A} 
\hat p_n(a|w)^{N_n(w,a)},
\end{equation}
where the empirical probabilities $\hat p_n(a|w)$ are given by 
\begin{equation}\label{tp}
\hat p_n(a|w) = \frac{N_n(w,a)}{N_n(w)}.
\end{equation}
Here and in the sequel we use the convention $0^0=1$, for example in
the case of $N_n(w,a) = 0$ in expression \ref{pmle}. Note that by 
Definition~\ref{feasible}, as $N_n(w)\geq 1$ for any $w\in\tau$, 
it is not necessary to give an extra definition
of $\hat p_n(a|w)$ in the case $N_n(w) = 0$.

Given a sequence $w$, with $N_n(w)\geq 1$, we will denote by
\begin{equation*}
\hat\P_{\ML,w}(x_1^n) =  \prod_{a\in A} \hat p_n(a|w)^{N_n(w,a)}.
\end{equation*}
Hence, we have 
\begin{equation*}
\hat\P_{\ML,\tau}(x_1^n) =  \prod_{w\in \tau}   \hat\P_{\ML,w}(x_1^n).
\end{equation*}

Let $f(n)$ be any positive function such that $f(n)\to +\infty$, when
$n\to +\infty$, and $n^{-1}f(n)\to 0$, when $n\to +\infty$. This
function will represent the generic penalizing term of our estimator,
replacing the function $\frac{|A|-1}{2}\log n$ in the classical definition of BIC \citep{csiszar2006}.
A function satisfying these conditions will be called \emph{penalizing term}.  

\begin{definition}
Given a penalizing term $f(n)$, the PL context
tree estimator is given by
\begin{equation}\label{est:bic}
\hat \tau(x_1^n) = \argmin{\tau\in\sft}\;\{\,-\log 
\hat\P_{\ML,\tau}(x_1^n) + |\tau|f(n)  \,\}.
\end{equation}
\end{definition}

As can be seen, the computation of the estimated context tree using its raw
definition would imply a search for the optimal tree on the set of all
feasible trees. This was the biggest drawback of this
approach, because the size of this set grows extremely fast as a
function of the maximal height $d$. Fortunately, there is a way of
computing the PL estimator without exploring the set of all trees, as
shown by \citet{csiszar2006}. The details of this algorithm are given in the Appendix and
will be used in the proof of our main result.
 
Let $K\in\N$.  Define the underestimation event 
with respect to the truncated tree $\tau_0|_{K}$ 
by 
\begin{equation*}
U_n^K = \bigcup_{w\in\Int(\tau_0|_K)}\{\,w\in \hat\tau_n(x_1^n)\}
\end{equation*}
and the overestimation event by
\begin{equation*}
O_n^K = \bigcup_{w\succ v\in\tau_0,\ell(v)<K}\;\{\, w\in\hat\tau_n(x_1^n)\}.
\end{equation*}

We are ready to present the main result in this paper. It
establishes upper bounds for the probability of occurrence of the underestimation and overestimation events.

\begin{theorem}\label{teo:bic}
Let $x_1,x_2,\dotsc$ be a sample of the stationary ergodic stochastic process $X_t$ having context tree $\tau_0$ and satisfying Assumption~1. For any constant $K\in\N$ there exist an integer $n_0$ and positive constants $c_1$, $c_2$, $c_3$ and $c_4$ depending on the process $X_t$ such that for any $n\geq n_0$
\begin{enumerate}
\item[\textup{(a)}] $\P\bigl[\,U_n^K]\;\leq \; c_1 \,e^{- c_2(n-d)}$;
\item[\textup{(b)}] $\P\bigl[\,O_n^K]\;\leq \; c_3 |A|^d\,e^{- c_4f(n)(\alpha_0^2/|A|)^d/d}$.
\end{enumerate}
\end{theorem}  
    
\begin{corollary}\label{cor:bic}
For any penalizing term $f(n)$ and any function $d(n)$ such that for any constant $c>0$,
\begin{equation}\label{serie}
\sum_{n\in\N} |A|^{d(n)}\exp[-\frac{f(n)c^{d(n)}}{d(n)}] < +\infty
\end{equation}
we have that there exists an integer $n_0$ depending on the process $X_t$ such that $\hat\tau_n(x_1^n)|_K=\tau_0|_K$ for any $n\geq n_0$. 
\end{corollary}  
  
\section{Proof of Theorem~\ref{teo:bic}}

Using Definition~\ref{def:max} and Lemma~\ref{prop:csis} and we see
that the tree in (\ref{est:bic}) can be written as
\begin{equation*}
  \hat\tau(x_1^n) = \{w\in \cup_{j=1}^dA^j\colon \X{w} = 0, \;
  \X{v} = 1\text{ for all } v\prec w \}
\end{equation*}
if $\X{\lambda} = 1$, and to $\{\lambda\}$ if $\X{\lambda} = 0$.  
Then, for $n$ sufficiently large in order to guarantee that  $\tau_0|_K$ will be in $\sft$ we have that 
\begin{equation*}
U_n^K = \bigcup_{w\in\Int(\tau_0|_K)} \{\,\X{w} = 0\,\} 
\end{equation*}
and
\begin{equation*}
O_n^K \subset \bigcup_{v\in\tau_0, \ell(v) < K}\;\{\, \X{v} = 1\,\}.
\end{equation*}
To prove (a)  let $w\in\Int(\tau_0|_K)$, then using Definition~\ref{def:xv} and Lemma~\ref{lem:csis} we have that
\begin{align*}
\P\bigl[\, \X{w} = 0 \,\bigr]  &\,= \,\P\bigl[\;\prod_{a\in A} 
\V{aw}\, \leq\, e^{-f(n)}\hat\P_{\ML,w}(x_1^n) \;\bigr]
\end{align*}
and for any $a\in A$
\[
\V{aw} = \max_{\tau\in\F^d_{aw}(x_1^n)} \prod_{s\in\tau}e^{-f(n)}\hat\P_{\ML,s}(x_1^n),
\]
where $\F^d_{aw}(x_1^n)$ is the set  containing all trees $\tau$ that have the form 
$\tau=\tau'\cap\{u\colon u\succeq aw\}$, with $\tau'\in\sft$.
Then 
\begin{align*}
\P\bigl[\, \X{w} = 0 \,\bigr]  \,= \,\P\bigl[\; \max_{\tau\in\F^d_{w}(x_1^n)} \prod_{s\in\tau}e^{-f(n)}\hat\P_{\ML,s}(x_1^n)
\, \leq\, e^{-f(n)}\hat\P_{\ML,w}(x_1^n) \;\bigr].
\end{align*}
For a tree $\tau\in\F^d_w(x_1^n)$ define the quantity
\begin{equation}\label{delta}
 \delta_{\tau}(w) =  \sum_{a\in A} \bigl[\sum_{u\in\tau}p(ua)\log p(a|u) -  
 p(wa)\log p(a|w)\bigr].
 \end{equation}
Using Jensen's inequality we can see that
$\delta_{\tau}(w) > 0$ unless $p(a|w)=p(a|u)$ for all $a\in A$ and all $u\in\tau$. 
Therefore, for a sufficiently large $n$ there must be a tree $\tau'_w\in\F^d_{w}(x_1^n)$ such that $\delta_{\tau'_w}(w) > 0$; if not 
we contradict the fact that $w\in\Int(\tau_0)$ and it is not a context in the sense of Definition~\ref{context}. Therefore 
\begin{align*}
\P\bigl[\, \X{w} = 0 \,\bigr]  \,\leq \,\P\bigl[\; \prod_{u\in\tau'_w}e^{-f(n)}\hat\P_{\ML,u}(x_1^n)
\, \leq\, e^{-f(n)}\hat\P_{\ML,w}(x_1^n) \;\bigr].
\end{align*}
Now we can apply the logarithm function on
both sides inside the probability obtaining that the right hand side equals  
\begin{equation*}
\P\bigl[\; \sum_{u\in\tau'_w}\log \hat\P_{\ML,u}(x_1^n) - 
\log\hat\P_{\ML,w}(x_1^n) \,\leq\, (|\tau'_w| - 1) f(n) \;\bigr].
\end{equation*}
Dividing by $n-d$
and  subtracting on both sides the term $\delta_{\tau'_w}(w)$
we have that for a  sufficiently large $n$ such that
\begin{equation*}\label{cond:n}
\frac{|\tau'_w|f(n)}{n-d} \, < \,\frac{\delta_{\tau'_w}(w)}{2}
\end{equation*}
 we can bound above the last expression by
\begin{align*}
\P\bigl[\, |L_n(w)|>\frac{\delta_{\tau'_w}(w)}{4}\bigr] \,+\,  \sum_{u\in\tau'_w}  \P\bigl[\, |L_n(u)|>\frac{\delta_{\tau'_w}(w)}{4}\bigr], 
\end{align*}
where for any finite sequence $s$
\[
L_n(s) =  \sum_{a\in A} p(sa)\log p(a|s) -  \frac{N_n(s,a)}{n-d} \log \hat p_n(a|s).
\]
Using Corollary~\ref{cor:estim1} we can bound above this expression by
\[
3\ee|A|^2(1+|\tau'_w|) \exp\bigl[-\frac{(n-d)\min(\delta_{\tau'_w},\delta_{\tau'_w}^2)\alpha_0^{2(h(\tau'_w)+1)}}
{1024e|A|^3(\alpha+\alpha_0)\log^2\alpha_0h(\tau'_w)}\bigr].
\]
We conclude the proof of part (a) by observing that we only have a finite number 
of sequences $w\in\Int(\tau_0|_K)$, so we can take
\[
c_1 = \max_{w\in\Int(\tau_0|_K)} \{3\ee|A|^2(1+|\tau'_w|)\}
\]
and
\[
c_2 = \min_{w\in\Int(\tau_0|_K)}\{\frac{\min(\delta_{\tau'_w},\delta_{\tau'_w}^2)\alpha_0^{2(h(\tau'_w)+1)}}{1024e|A|^3(\alpha+\alpha_0)\log^2\alpha_0h(\tau'_w)}\}.
\]
\noindent To prove part (b)  
observe that for any $w\in\tau_0$ with $\ell(w)<K$ 
\begin{equation}\label{x1}
\P\bigl[\, \X{w} = 1 \,\bigr]  = \P\bigl[\;\prod_{a\in A} \V{aw} \,>\, 
e^{-f(n)}\hat\P_{\ML,w}(x_1^n) \;\bigr].
\end{equation}
Using Lemma~\ref{lem:csis} we have that
\begin{equation*}
\prod_{a\in A} \V{aw}  =  \prod_{u\in \tau_w(x_1^n)}  e^{-f(n)}
\hat\P_{\ML,u}(x_1^n).
\end{equation*}
Then, applying the logarithm function the probability (\ref{x1})
is equal to
\begin{align}\label{exp:lem>}
&\P\bigl[\sum_{u\in \tau_w(x_1^n)} \log \, e^{-f(n)}
\hat\P_{\ML,u}(x_1^n)   >  \log e^{-f(n)}\hat\P_{\ML,w}(x_1^n) \;\bigr] \\
&= \,\P\bigl[\; \log\hat\P_{\ML,w}(x_1^n) - \!\!\sum_{u\in
\tau_w(x_1^n)}\!\!\log\hat\P_{\ML,u}(x_1^n) < (1
-|\tau_w(x_1^n)|) f(n) \;\bigr].\notag
\end{align}
We know, by the maximum likelihood estimator of the transition
probabilities that
\begin{equation}\label{mle}
\hat\P_{\ML,w}(x_1^n)\, \geq \, \prod_{a\in A} p(a|w)^{N_n(w,a)}.
\end{equation}
Therefore, we can bound above the right hand side of 
(\ref{exp:lem>}) by
\begin{align*}
\P\bigl[\,\sum_{a\in A} &N_n(w,a)\log p(a|w) -\!\!
\sum_{u\in\tau_w(x_1^n)}\!\! \log \hat\P_{\ML,u}(x_1^n)\,] <  
(1-|\tau_w(x_1^n)|) f(n) \bigr]\\
&=\,\P\bigl[\,\sum_{a\in A}\sum_{u\in\tau_w(x_1^n)}
N_n(u,a)\log\frac{p(a|u)}{\hat p_n(a|u)}\,] <
(1 - |\tau_w(x_1^n)|) f(n) \bigr].
\end{align*}
This equality follows by substituting $N_n(w,a)$ by $
\sum_{u\in\tau_w(x_1^n)} N_n(u,a)$ and the fact that $p(a|u) =
p(a|w)$ for all $u\in\tau_w(x_1^n)$, remembering that 
$w\in\tau_0$.  Observe that
\begin{align*}
\sum_{a\in A}\sum_{u\in\tau_w(x_1^n)} N_n(u,a)\log\frac{p(a|u)}
{\hat p_n(a|u)}  & \,
= \sum_{u\in\tau_w(x_1^n)} N_n(u)\sum_{a\in A} 
\hat p_n(a|u) \log\frac{p(a|u)}{\hat p_n(a|u)}\\
& \,= \, - \sum_{u\in\tau_w(x_1^n)} N_n(u)\, D\bigl(\hat
p_n(\cdot|u)\, \Vert\, p(\cdot|u)\bigr),
\end{align*}
where $D$ is the \emph{Kullback-Leibler divergence} between the
two distributions $\hat p_n(\cdot|u)$ and $p(\cdot|u)$ (see
the Appendix).  Using Lemma~\ref{lem:div} and dividing by $n-d$ we
have that
\begin{align*}
\P\bigl[\,-&\sum_{u\in\tau_w(x_1^n)}N_n(u)\, 
D\bigl(\hat p_n(\cdot|u)\, \Vert\, p(\cdot|u)\bigr)\,] < 
(1 - |\tau_w(x_1^n)|)f(n) \, \bigr]\\
&\, \leq\;\P\bigl[\,- \sum_{u\in\tau_w(x_1^n)}\frac{N_n(u)}{n-d}
\sum_{a\in A}\frac{[\hat p_n(a|u) - p(a|u)]^2}{p(a|u)}\,] < 
\frac{(1-|\tau_w(x_1^n)|)f(n)}{n-d}\, \bigr].
\end{align*}
As $\X{w}=1$ 
it follows that $|\tau_w(x_1^n)|>1$. On the other hand,
$N_n(u) \leq n-d$ and $f(n) > 0$. Therefore, we can bound
above the right hand side of the last expression by
\begin{equation*}
\sum_{u\in\tau_w(x_1^n)} \sum_{a\in A}\; \P\Bigl[\, \bigl|\hat p_n(a|u) - 
p(a|u)\bigr|\, > \sqrt{\frac{f(n)p(a|u)}{(n-d)|A||\tau_w(x_1^n)|}}\;
\Bigr].
\end{equation*}
Hence, using Corollary~\ref{cor:estim1} we can bound above this
expression by
\[
2\,\ee \, |A|^{d+2}\, \exp \bigl[\;-\, \frac{ f(n)\,\alpha_0^{2(d+1)}}{32e(\alpha+\alpha_0)|A|^{d+3}d}\bigl].
\]
This finishes the proof of Theorem~\ref{teo:bic}, by taking
\[
c_3 = 2\ee |A|^2\quad\text{ and }\quad c_4 = \frac{\alpha_0^2}{32e(\alpha+\alpha_0)|A|^3}.
\]

\noindent{\em Proof of Corollary~\ref{cor:bic}}.  It follows from the Borel-Cantelli Lemma and Theorem~\ref{teo:bic}, 
by noting that 
\[
\P[\hat\tau(x_1^n)|_K\neq \tau_0|_K] \leq \P[U_n^K]+\P[O_n^K]
\]
and the right hand side is summable in $n$ when condition (\ref{serie}) is satisfied. 
\section{Final Remarks}

The present paper presents upper bounds for the rate of convergence of 
penalized likelihood context tree estimators. 
We obtain an exponential bound for the underestimation event and 
 an under-exponential bound in the case of the overestimation event. 
These results generalizes the previous work by 
\cite{dorea2006}, who obtained similar bounds in the case of the estimation
of the order of a Markov chain, using also penalized likelihood criteria. 
One question that still remains open is if these bounds are optimal, as in the case
of an estimator introduced in \citet{finesso1996} for the estimation of the order of a Markov
chain. They prove that in the case of their estimator, the constant appearing in the 
underestimation bound is optimal, and that the overestimation bound can not be 
exponential if the estimator is universal, as in our case.  
The answer to these questions are important subjects for future work in this area.

\section{Appendix}

\subsection{The context tree maximizing principle}

The following definitions and results were taken from \citet{csiszar2006} and were included for completeness. 
 Definitions \ref{def:xv} and \ref{def:max} and Lemmas~\ref{lem:csis} and \ref{prop:csis} were originally proven for the usual penalizing term $f(n)=\frac{|A|-1}{2}\log n$, but can be adapted in a straightforward way to our setting.

Given two probability distributions $p$ and $q$ over $A$, the \emph{Kullback-Leibler divergence} 
is defined by
\begin{equation}\label{divergence}
D(p \Vert q) =   \sum_{a\in A} p(a) \log \frac{p(a)}{q(a)},
\end{equation}
where, by convention, $p(a) \log \frac{p(a)}{q(a)}$ equals $0$ if $p(a) = 0$ and $+\infty$
if $p(a) > q(a) = 0$. 

\begin{lemma}\label{lem:div}
If $p$ and $q$ are two probability distributions over $A$ then
\begin{equation}
D(p\Vert q) \,\leq \, \sum_{a\in A} \frac{[p(a) - q(a)]^2}{q(a)}.
\end{equation}
\end{lemma} 

\begin{proof}
See \citet[Lemma~6.3]{csiszar2006}.
\end{proof}

Consider the full tree $A^d$, and let $S^d$ denote the set of all
sequences of length at most $d$, that is $S^d=\cup_{j=0}^dA^j$.

\begin{definition}\label{def:xv}
  Given a sequence $w\in S^d$ with $N_n(w)\geq 1$, we define 
  recursively, starting from the sequences of the full tree $A^d$, the
  value
\begin{equation*}
  V_w(x_1^n) = \begin{cases}
    \max\{e^{-f(n)}\hat\P_{\ML,w}(x_1^n),\; \prod_{a\in A}V_{aw}(x_1^n)\}, 
    &\text{ if $0\leq \ell(w) < d,$}\\
    e^{-f(n)}\hat\P_{\ML,w}(x_1^n), & \text{ if $\ell(w) = d$}
\end{cases}
\end{equation*}
and the indicator
\begin{equation*}
  \X{w} = \begin{cases}
    1, &\text{ if $0 \leq  \ell(w) < d$ and }  
    \prod_{a\in A}V_{aw}(x_1^n) > e^{-f(n)}\hat\P_{\ML,w}(x_1^n),  \\
    0, & \text{ if $0 \leq \ell(w) < d$ and } 
    \prod_{a\in A}V_{aw}(x_1^n) \leq  e^{-f(n)}\hat\P_{\ML,w}(x_1^n), \\
    0, & \text{ if $\ell(w) = d$.}
\end{cases}
\end{equation*}
\end{definition}

\begin{definition}\label{def:max}
  Given $w\in S^d$ with $N_n(w)\geq 1$, the maximizing tree
  assign to the sequence $w$ is the tree
\begin{equation*}
  \tau_w(x_1^n) = \{u\in S^d\colon \X{u} = 0, \;\X{v} = 1
  \text{ for all }w\preceq v \prec u \}
\end{equation*}
if $\X{w} = 1$ and $\tau_w(x_1^n) = \{w\}$ if $\X{w} = 0$.
\end{definition}

For a sequence $w\in S^d$, with $N_n(w)\geq 1$, 
define $\F^d_w(x_1^n)$ as the set  containing all trees $\tau$ that have the form 
$\tau=\tau'\cap\{u\colon u\succeq w\}$, with $\tau'\in\sft$.

\begin{lemma}\label{lem:csis}
For any $w\in S^d$ with $N_n(w)\geq 1$,
\begin{equation*}
\V{w} =  \max_{\tau\in\F^d_w(x_1^n)}\prod_{u\in \tau}  e^{-f(n)}\hat\P_{\ML,u}(x_1^n) = \prod_{u\in \tau_w(x_1^n)}  e^{-f(n)}\hat\P_{\ML,u}(x_1^n).
\end{equation*} 
\end{lemma}
  
\begin{proof}
See \citet[Lemma~4.4]{csiszar2006}.
\end{proof}

\begin{lemma}\label{prop:csis}
  The context tree estimator $\hat\tau(x_1^n)$ in (\ref{est:bic})
  equals the maximizing tree assigned to the empty string $\lambda$,
  that is,
\begin{equation*}
\hat\tau(x_1^n)= \tau_\lambda(x_1^n).
\end{equation*}
\end{lemma} 

\begin{proof}
See \citet[Proposition~4.3]{csiszar2006}.
\end{proof}

From this result it follows that in order to obtain the tree maximizing the penalized maximum likelihood criteria  it is sufficient to assign to each sequence $w\in S^d$, with $N_n(w)\geq 1$, the indicator $\X{w}$ and then to get the maximizing tree $\tau_\lambda(x_1^n)$. The computational cost of this algorithm is linear in $n$ if $d(n)=o(n)$, as proven by \citet{csiszar2006}.   

\subsection{Exponential inequalities for empirical probabilities}

The following result was proven in \citet{galves2008}, we omit its proof here.  
\begin{theorem}\label{estim1} 
Assume the process $X_t$ satisfies Assumption~1, then for any 
finite sequence $w$, any symbol $a\in A$ and any $t>0$ the
following inequality holds
\begin{equation*}\label{Nn2}
\P(\,|N_n(w,a)-(n-d)p(wa)|\,>\,t\,)\,\leq \,e^{\frac1e} 
\exp \bigl[\frac{-t^2 C}{(n-d)\ell(wa)}\bigr],
\end{equation*}
where
\begin{equation*}\label{C}
C = \frac{\alpha_0}{8e(\alpha+\alpha_0)}.
\end{equation*}
\end{theorem}

As a consequence of Theorem~\ref{estim1} we obtain the
following corollary.

\begin{corollary}\label{cor:estim1} 
For any finite sequence $w$, with $p(w)>0$, any $t>0$ and any sufficiently large $n$ 
such that $N_n(w)\geq 1$ we have \\[-3mm]
\begin{itemize}
\item[\textup{(a)}] $\max_{a\in A} \P\bigl(|\hat{p}_n(a|w)-p(a|w)|>t \bigl)
\,\leq\,2\, \ee \, |A|\,\exp \bigl[-  
\frac{(n-d)t^2p(w)^2\alpha_0}{32e|A|^2(\alpha+\alpha_0)\ell(wa)}\bigl]\,;$\\
\item[\textup{(b)}]   $\P\bigl[\;| L_n(w) |> t\,\bigr] \,
\leq   3\, \ee\, |A|^2\, \exp\bigl[
-\frac{(n-d)\min(t,t^2)p(w)^2\alpha_0^2}{64e|A|^3
(\alpha+\alpha_0)\log^2\alpha_0\ell(wa)} \bigr]\,,$\\[2mm]
where $L_n(w) = \sum_{a\in A} p(wa)\log p(a|w) -  \frac{N_n(w,a)}{n-d} \log \hat p_n(a|w)$.
\end{itemize}
\end{corollary}

\begin{proof}
To prove (a) observe that
\[
p(a|w) = \frac{(n-d)p(wa)}{(n-d)p(w)}.
\]
Then, 
summing and substracting the term $\frac{N_n(w,a)}{(n-d)p(w)}$ we
obtain
\begin{align*}
\Bigl|\frac{N_n(w,a)}{N_n(w)} - \frac{(n-d)
p(wa)}{(n-d)p(w)}\Bigr| \,\leq\, 
&\frac{N_n(w,a)}{N_n(w)(n-d)p(w)}\,\bigl|(n-d)p(w) -
N_n(w)\bigr| \\
&\;+ \;\frac{1}{(n-d)p(w)}\;\bigl|N_n(w,a) - (n-d)p(wa)\bigr|.
\end{align*}
Therefore, as $\frac{N_n(w,a)}{N_n(w)} \leq 1$ we have
\begin{align*}
\P\bigl(|\hat{p}_n(a|w)-p(a|w)|>t \bigl) \; \leq&\; \P\Bigl(\bigl|(n-d)p(w) -
N_n(w)\bigr| > \frac{t(n-d)p(w)}{2}\Bigr)\\
& + \P\Bigl(\bigl|N_n(w,a) - (n-d)p(wa)\bigr| > \frac{t(n-d)p(w)}{2}\Bigr)
\end{align*}
We can write $N_n(w) = \sum_{b\in A} N_n(w,b)$ and $p(w) = \sum_{b\in A}
p(wb)$, then the right hand side of the last expression can be bounded above
by the sum
\begin{align*}
&\sum_{b\in A} \,\P\bigl(\,|N_n(w,b)-(n-d)p(wb)|\,>\,\frac{t (n-d)p(w)}{2|A|}
\bigl)\;\;+\\
&\qquad\qquad\qquad \P\bigl(\,|N_n(w,a)-(n-d)p(wa)|\,>\, \frac{t(n-d)
p(w)}{2} \bigl).
\end{align*}
Using Theorem~\ref{estim1} we can bound above this expression by
\begin{equation*}\label{s2}
\ee \,(|A|+1)\,\exp \bigl[- (n-d) \; \frac{t^2 p(w)^2 C}{4|A|^2
\ell(wa)}\bigl].
\end{equation*}
This finishes the proof of (a). 
To prove (b) observe that
\begin{align*}
\P\bigl[\;\bigl|L_n(w) \bigr| \,>\, t\,\bigr] \;\leq\;&\P\bigl[\;\bigl| \sum_{a\in A} \log p(a|w) \bigl(p(wa) -
\frac{N_n(w,a)}{n-d}\bigr)
\bigr| \,>\,\frac{t}{2}\,\bigr] \\
& + \P\bigl[\;\bigl| \sum_{a\in A} \frac{N_n(w,a)}{n-d} \log
\frac{p(a|w)}{\hat p_n(a|w)} \bigr| \,>\, \frac{t}{2}\,\bigr].
\end{align*}
Using Theorem~\ref{estim1} we have that
\begin{align}
\P\bigl[\;\bigl| \sum_{a\in A} \log p(a|w)&\bigl(p(wa) -
\frac{N_n(w,a)}{n-d}\bigr)
\bigr| \,>\,\frac{t}{2}\,\bigr] \notag\\
& \leq \; \sum_{a\in A} \,\P\bigl[\;\bigl| N_n(w,a) - (n-d)p(wa)
\bigr| \,>\,\frac{(n-d) t}{2\,|\!\log p(a|w)| |A|}\,\bigr]\notag
\\ \label{eq:1} & \leq \;\;\ee\,|A|\, \exp\bigl[
\frac{-(n-d)t^2C}{4|A|^2\log^2\alpha_0\ell(wa)} \bigr].
\end{align}
On the other hand, using the definition of the \emph{Kullback-Leibler
divergence}, Lemma~\ref{lem:div} and part (a) of this Corollary we obtain
\begin{align}
\P\bigl[\;\bigl| \sum_{a\in A} \frac{N_n(w,a)}{n-d}& \log
\frac{p(a|w)}{\hat p_n(a|w)} \bigr| \,>\, \frac{t}{2}\bigr] \,\leq
\, \P\bigl[\; D(\hat p(\cdot|w) || p(\cdot|w) )\,>\,
\frac{t}{2}\bigr]\notag \\
& \quad \leq\, \sum_{a\in A}\, \P\Bigl[\;\bigl| p(a|w) - \hat
p_n(a|w)
\bigr| \,>\,\sqrt{\frac{t p(a|w)}{2|A|} }\,\Bigr]\notag\\
& \quad \leq\, 2\, \ee\, |A|^2\, \exp\bigl[
-\frac{(n-d)tp(w)^2\alpha_0^2}{64e|A|^3
(\alpha+\alpha_0)\ell(wa)} \bigr].\label{eq:2}
\end{align}
Summing (\ref{eq:1}) and (\ref{eq:2}) we obtain the bound in part (b) and we conclude the proof of Corollary~\ref{cor:estim1}.
\end{proof}

\section*{Acknowledgments}

The author is thankful to Antonio Galves, Aur\'elien Garivier, Eric Moulines and 
Bernard Prum for interesting suggestions to improve the presentation of the results in this paper. 

\bibliography{./references}  

\bibliographystyle{dcu}

\end{document}